\documentclass[a4paper,12pt]{article}
\usepackage{graphicx}
\usepackage{a4wide}
\usepackage[english]{babel}
\usepackage{amsfonts}
\usepackage{amsmath}
\usepackage{amssymb}
\usepackage{dsfont}

\newtheorem{theorem}{Theorem}

\allowdisplaybreaks

\date{}

\DeclareMathOperator {\var}{Var_{[0,1]}}

\newcommand {\ve} {\varepsilon}

\numberwithin{equation}{section} \numberwithin{theorem}{section}
\numberwithin{lemma}{section} 
\numberwithin{corollary}{section}
\makeatletter\def\blfootnote{\xdef\@thefnmark{}\@footnotetext}\makeatother

\begin{document}
\title{\bf Lacunary sequences and permutations}
\author{C.\ Aistleitner\footnote{Graz University of Technology, Department for Analysis and Computational Number Theory, Steyrergasse 30, 8010 Graz, Austria \mbox{e-mail}: \texttt{aistleitner@math.tugraz.at}.
Research supported by FWF grant S9603-N23.}, I.\ Berkes\footnote{ Graz University of Technology, Institute of
Statistics, M\"unzgrabenstra{\ss}e 11, 8010 Graz, Austria.  \mbox{e-mail}: \texttt{berkes@tugraz.at}. Research
supported by the FWF Doctoral Program on Discrete Mathematics (FWF DK W1230-N13), FWF grant S9603-N23 and OTKA grants K 67961 and K 81928.} and R.\ Tichy\footnote{Graz University of Technology, Department for Analysis and Computational Number Theory, Steyrergasse 30, 8010 Graz, Austria. \mbox{e-mail}: \texttt{tichy@tugraz.at}.
Research supported by the FWF Doctoral Program on Discrete Mathematics (FWF DK W1230-N13) and FWF grant S9603-N23.}}

\maketitle \vskip0.5cm

\abstract{ By a classical principle of analysis, sufficiently thin
subsequences of general sequences of functions behave like
sequences of independent random variables. This observation not
only explains the remarkable properties of lacunary trigonometric
series, but also provides a powerful tool in many areas of
analysis. In contrast to ``true'' random processes, however, the
probabilistic structure of lacunary sequences is  not
permutation-invariant and the analytic properties of such
sequences can change radically after rearrangement. The purpose of
this paper is to survey some recent results of the authors on
permuted function series. We will  see that rearrangement
properties of lacunary trigonometric series $\sum (a_k\cos
n_kx+b_k \sin n_kx)$ and their nonharmonic analogues $\sum c_k
f(n_kx)$ are intimately connected with the number theoretic
properties of $(n_k)_{k \geq 1}$  and we will give a complete
characterization of permutational invariance in terms of the
Diophantine properties of $(n_k)_{k \geq 1}$. We will also see
that in a certain statistical sense, permutational invariance is
the ``typical'' behavior of lacunary sequences.}



\section{Introduction}

Let $(n_k)_{k \geq 1}$ be a sequence of positive integers
satisfying the Hadamard gap condition
\begin{equation}\label{had}
n_{k+1}/n_k\ge q>1 \qquad (k=1, 2, \ldots).
\end{equation}
Salem and Zygmund \cite{sz1947} proved
that if $(a_k)_{k \geq 1}$ is a sequence of real numbers
satisfying
\begin{equation}\label{an}
a_N=o(A_N) \quad  \text{with} \quad
A_N=\frac{1}{2}\left(\sum_{k=1}^N a_k^2\right)^{1/2},
\end{equation}
then $(\cos 2\pi n_kx)_{k\ge 1}$ obeys the central limit theorem
\begin{equation}\label{cosclt}
\lim_{N\to\infty} \lambda \{x\in (0, 1):A_N^{-1} \sum_{k=1}^N
a_k\cos 2\pi n_kx\le t\}=(2\pi)^{-1/2}  \int_{-\infty}^t
e^{-u^2/2}du,
\end{equation}
where $\lambda$ denotes the Lebesgue measure. Under the same gap
condition Weiss \cite{wetrig} proved (cf.\ also Salem and Zygmund
\cite{sz1950}, Erd\H{o}s and G\'al \cite{eg}) that if $(a_k)_{k
\geq 1}$ satisfies
\begin{equation}\label{kol}
a_N=o(A_N/(\log\log A_N)^{1/2})
\end{equation}
then $(\cos 2\pi n_kx)_{k \geq 1}$ obeys the law of the iterated
logarithm
\begin{equation}\label{coslil}
\limsup_{N\to\infty} \, (2A_N^2\log\log A_N)^{-1/2} \sum_{k=1}^N
a_k\cos 2\pi n_kx=1 \qquad \textup{a.e.}
\end{equation}
Comparing these results with the classical forms of the central
limit theorem and law of the iterated logarithm in probability
theory, we see that under the gap condition (\ref{had}) the
functions $\cos 2\pi n_kx$ behave like independent random
variables. Using martingale techniques, Philipp and Stout
\cite{phs} proved that if instead of (\ref{an}) we assume
$a_N=o(A_N^{1-\delta})$ for some $\delta>0$, then on the
probability space $([0, 1], {\cal B}, \lambda)$ there exists a
Brownian motion process $\{W(t), \, t\ge 0\}$ such that
\begin{equation}\label{br}
\sum_{k=1}^{N} \cos 2\pi
n_kx=W(A_N)+O\left(A_N^{1/2-\varrho}\right)\qquad \text{a.s.}
\end{equation}
for some $\varrho > 0$. The last relation implies not only the CLT
and LIL for $(\cos 2\pi n_kx)_{k \geq 1}$, but a whole class of
further limit theorems for  independent random variables; for
examples and discussion we refer to \cite{phs}.

The previous results extend, in a modified form, to lacunary
subsequences of the system $\{f(nx)\}_{n\ge 1}$ where $f$ is a
periodic measurable function, but the asymptotic properties of
this system are much more complicated than those of the
trigonometric system. By a conjecture of Khinchin \cite{kh},  if
$f$ has period 1 and is Lebesgue integrable on $(0,1)$, then
\begin{equation}\label{khi}
\lim_{N\to \infty} \frac{1}{N}\sum_{k=1}^N f(kx)=\int_0^1 f(t)\,
dt \quad \textup{a.e.}
\end{equation}
This remained open for almost 50 years until Marstrand \cite{mar}
disproved it, but even today, no precise condition for the
validity of  (\ref{khi}) is known. Similarly, there is no analogue
of Carleson's theorem \cite{car} for the system $(f(nx))_{n \geq
1}$ and we do not know under what conditions the series
$\sum_{k=1}^\infty c_k f(kx)$ converges almost everywhere. In the
lacunary case, Kac \cite{ka} proved that if $f$ satisfies a
Lipschitz condition, then $f(2^kx)$ obeys  a central limit theorem
similar to (\ref{cosclt}) and not much later, Erd\H{o}s and Fortet
(see \cite{kac1949}, p.\ 655) showed that the CLT fails for
$f(n_kx)$ for  $n_k=2^k-1$ even for some trigonometric polynomials
$f$. Gaposhkin \cite{gap1966} proved that $f(n_kx)$ obeys the CLT
if $n_{k+1}/n_k\to\alpha$ where $\alpha^r$ is irrational for $r=1,
2 \ldots$ and the same holds if all the fractions $n_{k+1}/n_k$
are integers. He also showed (see \cite{gap1970}) that the
validity of the CLT for $f(n_kx)$ is closely related to the number
of solutions of the Diophantine equation
\begin{equation}\label{diogap}
an_k+bn_\ell=c, \qquad 1\le k, \ell \le N.
\end{equation}
Improving these results, Aistleitner and Berkes \cite{aibe}
recently gave a necessary and sufficient Diophantine condition for
the CLT for $f(n_kx)$. As the proofs of these results show,  the
asymptotic behavior of $f(n_kx)$ is determined by a complicated
interplay between  the arithmetic properties of $(n_k)_{k \geq 1}$
and the Fourier coefficients of $f$ and the combination of
probabilistic and number-theoretic effects leads to a unique,
highly interesting asymptotic behavior. Let
$$
D_N (x_1, \dots, x_N): = \sup_{0 \leq a < b < 1} \left|
\frac{\sum_{k=1}^N \mathds{1}_{[a,b)} (x_k)}{N} - (b-a) \right|
$$
denote the discrepancy (mod 1) of the finite sequence $(x_1,
\ldots, x_N)$, where $\mathds{1}_{[a,b)}$  is the indicator
function of the interval $[a,b)$, extended periodically to
$\mathbb R$. Philipp \cite{plt} proved that if $(n_k)_{k \geq 1}$
satisfies the Hadamard gap condition (\ref{had}), then the
discrepancy $D_N(n_kx)$ of the sequence $\{n_kx, 1\le k \le N\}$
obeys the LIL
\begin{equation}\label {phlil}
 \frac{1}{4\sqrt{2}} \leq \limsup_{N \to \infty} \frac{N D_N(n_kx)}
  {\sqrt{2N \log \log N}}
 \leq C_q \quad \textup{a.e.},
\end{equation}
where $C_q$ is a number depending on $q$. Note that if $(\xi_k)_{k
\geq 1}$ is a sequence of independent random variables with
uniform distribution over $(0, 1)$, then
\begin{equation}\label{cslil}
 \limsup_{N\to\infty}  \frac{ND_N(\xi_k)}{\sqrt{2N\log\log N}}
 =\frac{1}{2} \qquad
\end{equation}
with probability one by the Chung-Smirnov LIL (see e.g.\
\cite{sw}, p.\ 504). A comparison of (\ref{phlil}) and
(\ref{cslil}) shows again that the sequence $(n_kx)_{k \geq 1}$
mod 1 behaves like a sequence of i.i.d.\  random variables.
Surprisingly, however, the limsup in (\ref{phlil}) can be
different from the constant $1/2$ in (\ref{cslil}) and, as
Fukuyama \cite{ft} showed, it depends sensitively on $(n_k)_{k
\geq 1}$. For example, for $n_k=a^k$, $a\ge 2$ the limsup
$\Sigma_a$  in (\ref{phlil}) equals

\bigskip
\begin{eqnarray}
\Sigma_a & = & \sqrt{42}/9  \phantom{999999999999} \textrm{if} ~
a =2 \nonumber\\
\Sigma_a & = & \frac{\sqrt{(a+1)a(a-2)}}{2 \sqrt{(a-1)^3}} \quad
\textrm{if} ~ a \geq 4 ~
\textrm{is an even integer},\nonumber\\
\Sigma_a & = & \frac{\sqrt{a+1}}{2 \sqrt{a-1}}
\phantom{9999999999} \textrm{if} ~ a \geq 3 ~ \textrm{is an odd
integer} \nonumber.
\end{eqnarray}
It is even more surprising that, as Fukuyama \cite{ft2} showed,
the limsup in (\ref{phlil}) is not permutation-invariant and can
change after a rearrangement of $(n_k)_{k \geq 1}$. Similarly,
$$\limsup_{N\to\infty}\, (N\log\log N)^{-1/2} \sum_{k=1}^N
f(n_kx)$$ and the limiting variance in the CLT for $N^{-1/2}
\sum_{k=1}^N f(n_kx)$ can change if we permute the sequence
$(n_k)_{k \geq 1}$. These results show that even though lacunary
subsequences of $(f(nx))_{n \geq 1}$ satisfy a large class of
limit theorems for i.i.d.\ random variables  and an i.i.d.\
sequence is a symmetric structure,  the behavior of lacunary
sequences is generally nonsymmetric. The purpose of the present
paper is to give a detailed analysis of the probabilistic
structure of $f(n_kx)$ and to clear up the effect of permutations
on its asymptotic properties. The proofs of our results will be
given in \cite{aibeti1},  \cite{aibeti2}, \cite{aibeti3}.

\section{The trigonometric case}

By Carleson's theorem \cite{car}, if $f\in L_2(0, 2\pi)$ then its
Fourier series  \begin{equation}\label {fou} f\sim
\frac{a_0}{2}+\sum_{k=1}^\infty \left(a_k\cos kx +b_k\sin
kx\right)
\end{equation}
converges almost everywhere. However, as was shown by Kolmogorov
(see \cite{kome}), there exists an $f\in L^2(0, 2\pi)$ whose
Fourier series (\ref{fou}) diverges a.e.\ after a suitable
permutation of its terms. This shows that the asymptotic
properties of the trigonometric system $\{\cos kx, \sin kx\}_{k\ge
1}$ are not permutation-invariant. On the other hand, Zygmund
\cite{z30} proved that if $(n_k)_{k \geq 1}$ satisfies the
Hadamard gap condition (\ref{had}) and
\begin{equation}\label{ell2}
\sum_{k=1}^\infty \left(a_k^2+b_k^2\right)<\infty
\end{equation}
then
\begin{equation}\label {lactr}
\sum_{k=1}^\infty (a_k\cos n_kx +b_k\sin n_kx)
\end{equation}
converges almost everywhere after any rearrangement of its terms,
giving a per\-mu\-tat\-ion-invariant property of lacunary
trigonometric series. Our first result below states that under
(\ref{had}) the systems $(\cos n_kx)_{k \geq 1}$, $(\sin n_kx)_{k
\geq 1}$ satisfy also the central limit theorem and law of the
iterated logarithm in a permutation-invariant form. More
precisely, we have

\begin{theorem}\label{permcltlil}
Let $(n_k)_{k \geq 1}$ be a sequence of positive integers
satisfying (\ref{had}) and let $\sigma: {\mathbb N}\to{\mathbb N}$
be a permutation of the positive integers. Then we have
\begin{equation}\label{permclta}
\lim_{N\to\infty} \lambda \{x\in (0, 1): \sum_{k=1}^N \cos 2\pi
n_{\sigma(k)}x\le t\sqrt{N/2}\}=(2\pi)^{-1/2} \int_{-\infty}^t
e^{-u^2/2}du
\end{equation}
and
\begin{equation}\label{permlil}
\limsup_{N\to\infty} \, (N\log\log N)^{-1/2}\sum_{k=1}^N \cos 2\pi
n_{\sigma(k)}x=1 \qquad \textup{a.e.}
\end{equation}
\end{theorem}

Note that for the unpermuted CLT and LIL we need much weaker gap
conditions than (\ref{had}). In fact, Takahashi \cite{tak68b},
\cite{tak72} (cf.\ also Erd\H{o}s \cite{er62}) showed that if a
sequence $(n_k)_{k \geq 1}$ of integers satisfies
\begin{equation}\label{ergap}
n_{k+1}/n_k\ge 1+k^{-\alpha}, \qquad 0\le \alpha<1/2
\end{equation}
then for any sequence $(a_k)_{k \geq 1}$ satisfying
$$
a_N=o(A_N N^{-\alpha}) \quad  \text{with} \quad
A_N=\frac{1}{2}\left(\sum_{k=1}^N a_k^2\right)^{1/2}
$$
we have the CLT (\ref{cosclt}) and LIL (\ref{coslil}). Note,
however, that (\ref{ergap}) does not imply permutation-invariance
and the following result shows that permutation-invariance fails
under any gap condition weaker than (\ref{had}).

\begin{theorem}\label{counterex}
For any positive sequence $(\varepsilon_k)_{k \geq 1}$ tending to
0, there exists a sequence $(n_k)_{k \geq 1}$ of positive integers
satisfying
$$n_{k+1}/n_k\ge 1+\varepsilon_k, \qquad k\ge k_0$$
and a permutation $\sigma: {\mathbb N}\to {\mathbb N}$ of the
positive integers such that the permuted central limit theorem
(\ref{permclta}) and the permuted law of the iterated logarithm
(\ref{permlil}) fail.
\end{theorem}

By a theorem of Erd\H{o}s \cite{er1941}, if $(n_k)_{k \geq 1}$ is
any (not necessarily increasing) sequence of different positive
integers such that for any integer $\nu>0$ the number of solutions
of the Diophantine equation
$$n_k \pm n_\ell=\nu, \qquad k, \ell \ge 1$$
is bounded by a constant $C$ independent of $\nu$, then the series
(\ref{lactr}) converges a.e.\ provided (\ref{ell2}) holds.  Since
this Diophantine property is permutation invariant, it implies the
a.e.\ unconditional convergence of (\ref{lactr}) as well.  Note
that Erd\H{o}s' condition  is much weaker than (\ref{had}); in
fact, it holds even for some polynomially growing sequences
$(n_k)_{k \geq 1}$. How slowly a sequence $(n_k)_{k \geq 1}$
satisfying this condition can grow is a well known open problem in
number theory; see Halberstam and Roth \cite{haro}, p.\ 234 and
Ajtai et al. \cite{aj}.

\section{The system $f(nx)$}

Let $f$ be a measurable function satisfying
\begin{equation} \label{fcond1}
f(x+1)=f(x), \quad \int_0^1 f(x)~dx=0, \quad \int_0^2 f^2(x)\,
dx<\infty
\end{equation}
and let $(n_k)_{k \geq 1}$ be a sequence of integers satisfying
the Hadamard gap condition (\ref{had}). The central limit theorem
for $f(n_kx)$ has a long history discussed in Section 1. To
formulate criteria for the permutation-invariant CLT and LIL, let
us say that a sequence $(n_k)_{k \geq 1}$ of positive integers
satisfies

\bigskip\noindent
{\bf Condition ${\mathbf B}_2$}, if for any fixed  nonzero
integers $a, b, c$ the number of solutions $(k, l)$ of the
Diophantine equation
\begin{equation}\label{gap}
an_k+bn_l=c
\end{equation}
is bounded by a constant $K(a, b)$, independent of $c$.

\bigskip\noindent
{\bf Condition ${\mathbf B}^{(s)}_2$} (strong $\mathbf{B}_2$), if
for any fixed integers $a\ne 0$, $b\ne 0$, $c$ the number of
solutions $(k, l)$ of the Diophantine equation (\ref{gap}) is
bounded by a constant $K(a, b)$, independent of $c$, where for
$c=0$ we require also $k\ne l$.

\bigskip

\noindent {\bf Condition ${\mathbf B}_2^{(w)}$} (weak ${\bf
B}_2$), if for any fixed nonzero integers $a, b, c$ the number of
solutions $(k, l)$ of the Diophantine equation
\begin{equation}\label{gap2}
an_k+bn_l=c, \qquad 1\le k, l\le N
\end{equation}
is $o(N)$, uniformly in $c$. \\

\bigskip

\noindent {\bf Condition ${\mathbf B}_2^{(0)}$}, if for any fixed
nonzero integers $a, b$ the number of solutions $(k, l)$ of the
Diophantine equation
\begin{equation}\label{gap0}
an_k+bn_l=0, \qquad 1\le k, l\le N, \ k\ne l
\end{equation}
is $o(N)$. \\

Condition ${\mathbf B}_2$ was introduced by Sidon \cite{si} in his
investigations of trigonometric series.  Gaposhkin \cite{gap1970}
proved that under mild smoothness assumptions on $f$, condition
${\mathbf B}_2$ implies the CLT for $f(n_kx)$ and Berkes and
Philipp \cite{beph2} showed that the same condition also implies a
Wiener approximation for the partial sums of $f(n_kx)$, similar to
(\ref{br}). Recently,  Aistleitner and Berkes \cite{aibe} proved
that the CLT holds for $f(n_kx)$ also under ${\mathbf B}_2^{(w)}$
and this condition is necessary. This settles the CLT problem for
$f(n_kx)$, but, as we noted, the validity of the CLT does not
imply permutation-invariant behavior of $f(n_kx)$. The purpose of
this section is to give a precise description of the CLT and LIL
behavior of permuted sums $\sum_{k=1}^N f(n_{\sigma(k)}x)$ and in
particular, to obtain characterizations of permutation invariance.

Our first result shows that if we assume the slightly stronger gap
condition
\begin{equation}\label{infgap}
n_{k+1}/n_k\to\infty
\end{equation}
then the behavior of $f(n_kx)$ is permutation-invariant,
regardless the number theoretic structure of $(n_k)_{k \geq 1}$.
In what follows, let $\|\cdot \|$ denote the $L_2(0, 1)$ norm.

\begin{theorem}\label{th1}
Let $(n_k)_{k \geq 1}$ be a sequence of positive integers
satisfying the gap condition (\ref{infgap}). Then for any
permutation $\sigma: {\mathbb N}\to {\mathbb N}$ of the integers
and for any measurable function $f:{\mathbb R}\to{\mathbb R}$
satisfying
\begin{equation} \label{fcond}
f(x+1)=f(x), \quad \int_0^1 f(x)~dx=0, \quad \var f < + \infty
\end{equation}
we have
\begin{equation}\label{permcltx}
{\frac{1}{\sqrt{N}}}\sum_{k=1}^N f(n_{\sigma (k)} x)
\stackrel{\mathcal D}
{\rightarrow} \mathcal{N}(0, \|f\|^2)\\
\end{equation}
and
\begin{equation}\label{lilperm1}
\limsup_{N\to \infty} \frac{1}{\sqrt{2N\log\log N}} \sum_{k=1}^N
f(n_{\sigma (k)} x)=\|f\| \quad\textup{a.e.}\\
\end{equation}
Moreover, for any permutation $\sigma$ of $\mathbb{N}$ we have
\begin{equation}\label {lilinfgap} \limsup_{N \to \infty}
\frac{ND_N(n_{\sigma(k)}x)}{\sqrt{2 N \log \log N}} =\frac{1}{2}
\quad \textup{a.e.}
\end{equation}
\end{theorem}

Our next theorem shows that if we slightly strengthen
(\ref{infgap}) then not only the CLT and LIL, but a much larger
class of limit theorems becomes permutation-invariant.

\begin{theorem}\label{th2} Let $f$ be a function satisfying (\ref{fcond1}) and
the Lipschitz $\alpha$ condition. Let $(n_k)_{k \geq 1}$ be an
increasing sequence of positive numbers such that
\begin{equation}\label{nka}
\sum_{k=1}^\infty (n_k/n_{k+1})^\alpha<\infty.
\end{equation}
Then there exists a bounded i.i.d.\ sequence $(g_k)$ of functions
on $(0, 1)$ such that
\begin{equation}\label{almiid}
\sum_{k=1}^\infty |f(n_kx)-g_k(x)|<\infty \qquad \textup{a.e.}
\end{equation}
\end{theorem}

Let $\sigma$ be a permutation of $\mathbb N$. Relation
(\ref{almiid}) implies that
$$\sum_{k=1}^\infty \left|f(n_{\sigma(k)}x)-g_{\sigma(k)}(x)\right|<\infty
\qquad \textup{a.e.}
$$
and consequently
\begin{equation}\label{equiv}
\sum_{k=1}^N f(n_{\sigma(k)}x)-\sum_{k=1}^N g_{\sigma(k)}(x) =O(1)
\qquad \textup{a.e.}
\end{equation}
Since the i.i.d.\ sequences $(g_k)$ and $(g_{\sigma(k)})$ are
probabilistically equivalent, relation (\ref{equiv}) implies that,
up to an error term $O(1)$, the asymptotic properties of the
partial sums $\sum_{k=1}^N f(n_{\sigma(k)}x)$ are the same. Thus
Theorem \ref{th2} expresses a very strong form of permutation
invariance of the sequence $f(n_kx)$. Condition (\ref{nka}) is
satisfied e.g.\ if $n_k= 2^{[ck \log_2 k]}$ with $c>1/\alpha$.

The proof of Theorem \ref{th2} shows that the approximating
i.i.d.\ sequence $(g_k)$ can be chosen to satisfy
\begin{equation}\label{inprob}
\mu \{x\in (0, 1): |f(n_kx)-g_k(x)| \ge \ve_k\}\le \ve_k, \qquad
k=1,2, \ldots
\end{equation}
with $\ve_k=(n_{k}/n_{k+1})^\alpha$. This gives more precise
information than (\ref{almiid}) if $(n_k)_{k \geq 1}$ grows very
rapidly. Actually, the approximation given by (\ref{inprob}) is
best possible. Let $(n_k)_{k \geq 1}$ be an increasing sequence of
positive integers such that the ratios $n_{k+1}/n_k$ are integers
and $\sum_{k=1}^\infty (n_k/n_{k+1})=\infty$. Then there exists no
i.i.d.\ sequence $(g_n)$ of functions on $[0, 1]$  such that
\begin{equation}\label{trigiid}
\mu \{x: |\cos 2\pi n_kx-g_k(x)| \ge \ve_k\}\le \ve_k, \qquad k=1,
2, \ldots.
\end{equation}
with $\sum_{k=1}^\infty  \ve_k<\infty$.\\

So far, we investigated the permutational invariance of $f(n_kx)$
under the growth condition $n_{k+1}/n_k\to\infty$. Assuming only
the Hadamard gap condition (\ref{had}), the situation becomes more
complex and the number theoretic structure of $(n_k)_{k \geq 1}$
comes into play. Our first result gives a necessary and sufficient
condition for the permuted partial sums $\sum_{k=1}^N
f(n_{\sigma(k)}x)$ to have only Gaussian limit distributions and
gives precise criteria this to happen for a specific permutation
$\sigma$. \\

\begin{theorem}\label{th3}
Let $(n_k)_{k \geq 1}$ be a sequence of positive integers
satisfying the Hadamard gap condition (\ref{had}) and condition
$\mathbf{B}_2$. Let $f$ satisfy (\ref{fcond}) and let $\sigma$ be
a permutation of $\mathbb N$. Then $ N^{-1/2}\sum_{k=1}^N
f(n_{\sigma (k)}x)$ has a limit distribution iff
\begin{equation}\label{gamma}
\gamma= \lim_{N\to\infty} N^{-1}\int_0^1 \left(\sum_{k=1}^N
f(n_{\sigma (k)}x)\right)^2 dx \
\end{equation}
exists, and then
\begin{equation}\label{fclt}
 N^{-1/2}\sum_{k=1}^N f(n_{\sigma (k)}x)\to_d N(0, \gamma).
\end{equation}
(If $\gamma=0$ then the limit distribution is degenerate.)
\end{theorem}

Theorem \ref{th3} is best possible in the following sense:
\begin{theorem} \label{th3bp}
If condition $\mathbf{B}_2$ fails, there exists a permutation
$\sigma$ of $\mathbb N$ such that the limit in (\ref{gamma})
exists, but the normed partial sums in (\ref{fclt}) do not have a
Gaussian limit distribution.
\end{theorem}

In other words, under the Hadamard gap condition and condition
$\mathbf{B}_2$,  the limit distribution of $N^{-1/2}\sum_{k=1}^N
f(n_{\sigma (k)} x)$ can only be Gaussian, but the variance of the
limit distribution depends on the constant $\gamma$ in
(\ref{gamma}) which, as simple examples show, is not
permutation-invariant. For example, if $n_k=2^k$ and $\sigma$ is
the identity  permutation, then (\ref{gamma}) holds with
\begin{equation}\label{kac}
\gamma=\gamma_f=\int_0^1 f^2(x)dx +2\sum_{k=1}^\infty \int_0^1
f(x)f(2^kx)dx
\end{equation}
(see Kac \cite{ka}). Using an idea of Fukuyama \cite{ft2}, one can
construct permutations $\sigma$ of $\mathbb N$ such that
\begin{equation}
\lim_{N\to\infty} \frac{1}{N} \int_0^1 \left(\sum_{k=1}^N
f(n_{\sigma(k)} x)\right)^2 dx =\gamma_{\sigma, f}
\end{equation}
with $\gamma_{\sigma, f} \ne \gamma_f$. Actually, the set of
possible values $\gamma_{\sigma, f}$ belonging to all permutations
$\sigma$ contains the interval $I_f=[\gamma_f, \|f\|^2]$ and it is
equal to this interval provided the Fourier coefficients of $f$
are nonnegative. For general $f$ this is
false (for details, see Aistleitner, Berkes and Tichy \cite{abt}).\\

Under the slightly stronger condition $\mathbf{B}_2^{(s)}$ we have
\begin{theorem}\label{co1} Let $(n_k)_{k \geq 1}$ be a sequence of positive
integers satisfying the Hadamard gap condition (\ref{had}) and
condition $\mathbf{B}_2^{(s)}$. Let $f$ satisfy (\ref{fcond}) and
let $\sigma$ be a  permutation of $\mathbb N$. Then the central
limit theorem (\ref{fclt}) holds with $\gamma=\|f\|^2$.
\end{theorem}

We now pass to the problem of the LIL.
\begin{theorem} \label{thlil}
Let $(n_k)_{k \geq 1}$ be a sequence of positive integers
satisfying the Hadamard gap condition (\ref{had}) and condition
$\mathbf{B}_2$. Let $f$ be a measurable function satisfying
(\ref{fcond}), let $\sigma$ be a permutation of $\mathbb N$ and
assume that the limit
(\ref{gamma}) exists. Then we have
\begin{equation}\label{lilperm3}
\limsup_{N\to \infty} \frac{\sum_{k=1}^N f(n_{\sigma (k)}
x)}{\sqrt{2N\log\log N}} =\gamma ^{1/2}\quad\textup{a.e.}\\
\end{equation}
\end{theorem}

\begin{theorem} \label{thlils}
Let $(n_k)_{k \geq 1}$ be a sequence of positive integers
satisfying  the Hadamard gap condition (\ref{had}) and condition
$\mathbf{B}_2^{(s)}$. The for any measurable function satisfying
(\ref{fcond}) and any permutation $\sigma$ of $\mathbb N$ we have
$$
\limsup_{N\to \infty} \frac{\sum_{k=1}^N f(n_{\sigma (k)}
x)}{\sqrt{2N\log\log N}} = \|f\| \quad\textup{a.e.}\\
$$
\end{theorem}

The proof of Theorems \ref{th3} and \ref{thlil} shows that if $f$
is a trigonometric polynomial of degree $d$, then in conditions
$\mathbf{B}_2$ resp.\ $\mathbf{B}_2^{(s)}$ it suffices to have the
bound for the number of solutions of (\ref{gap}) for coefficients
$a,b$ satisfying $|a|\leq d, |b| \leq d$. Applying this with $d=1$
and using the the fact that for a Hadamard lacunary sequence
$(n_k)_{k \geq 1}$ and $c \in \mathbb{Z}$ the number of solutions
$(k,l)$, $k\ne l$  of
$$
n_k \pm n_l = c
$$
is bounded by a constant which is independent of $c$ (see Zygmund
\cite[p.\ 203]{zt}), we get Theorem \ref{permcltlil} of the
previous section.\\

Theorem \ref{th3bp} shows that condition $\mathbf{B}_2$ is best
possible in Theorem \ref{th3}. We were not able to decide whether
this condition is also best possible Theorem \ref{thlil}, but
condition $\mathbf{B}_2$ is nearly best possible in Theorem
\ref{thlil} in the following sense: if there exist nonzero
integers $a,b,c$ such that the Diophantine equation
$$
a n_k + b n_l = c
$$
has infinitely many solutions $(k,l)$ with $k \neq l$, then the
LIL for $f(n_{\sigma(k)} x)$ fails to hold for a suitable
permutation $\sigma$ and a suitable trigonometric polynomial $f$.

\begin{theorem}\label{dlil}
Let $(n_k)_{k \geq 1}$ be a sequence of positive integers
satisfying (\ref{had}) and condition $\mathbf{B}_2^{(s)}$. Then
for any permutation $\sigma: {\mathbb N}\to {\mathbb N}$ we have
\begin{equation}\label {lilinfgap2} \limsup_{N \to \infty}
\frac{ND_N(n_{\sigma(k)}x)}{\sqrt{2 N \log \log N}} =\frac{1}{2}
\quad \textup{a.e.}\\
\end{equation}
\end{theorem}

All the results formulated so far assumed the Hadamard gap
condition (\ref{had}) or the stronger condition (\ref{infgap}). If
we weaken (\ref{had}), i.e. we allow subexponential sequences
$(n_k)_{k \geq 1}$, we need much stronger Diophantine conditions
even for the unpermuted CLT and LIL for $f(n_kx)$. Specifically,
we need uniform bounds for the number of solutions of Diophantine
equations of the type
\begin{equation}\label{dio}
a_1 n_{k_1}+\ldots +a_p n_{k_p}=b.
\end{equation}
Call a solution of (\ref{dio}) {\it nondegenerate} if no subsum of
the left hand side equals 0. Let us say that a sequence $(n_k)_{k
\geq 1}$ of positive integers satisfies

\bigskip\noindent
{\bf Condition ${\mathbf A}_p$}, if there exists a constant
$C_p\ge 1$  such that for any integer $b\ne 0$ and any nonzero
integers $a_1, \ldots, a_p$  the number of nondegenerate solutions
of the
Diophantine equation (\ref{dio}) is at most $C_p$. \\

The following results are the analogues of Theorems
\ref{th3}--\ref{dlil}  without growth conditions on $(n_k)_{k \geq
1}$.

\begin{theorem}\label{th4}
Let $(n_k)_{k \geq 1}$ be an increasing sequence of positive
integers satisfying condition $\mathbf{A}_p$ for all $p\ge 2$. Let
$f$ satisfy (\ref{fcond}), let $\sigma$ be a permutation of
$\mathbb N$ and assume that the limit (\ref{gamma}) exists. Then
the permuted CLT (\ref{fclt}) is valid.
\end{theorem}

\begin{theorem}\label{th5}
Let $(n_k)_{k \geq 1}$ be an increasing sequence of positive
integers satisfying condition ${\mathbf A}_p$ for all $p\ge 2$
with $C_p \le \exp(Cp^\alpha)$ for some $\alpha > 0$. Moreover,
assume that $f$ satisfies (\ref{fcond}), $\sigma$ is a permutation
of $\mathbb N$ and (\ref{gamma}) holds. Then the permuted LIL
(\ref{lilperm3}) is valid.
\end{theorem}

Note that for the validity of the LIL we require a specific bound
for the constants $C_p$ in condition $\mathbf{A}_p$. For
subexponentially growing $(n_k)_{k \geq 1}$, verifying property
$\mathbf{A}_p$ is a difficult number-theoretic problem. Classical
examples of such sequences are the Hardy-Littlewood-P\'olya
sequences, i.e.\ increasing sequences $(n_k)_{k \geq 1}$
consisting of all positive integers of the form
$q_1^{\alpha_1}\cdots q_\tau^{\alpha_\tau}$ ($\alpha_1, \ldots
\alpha_\tau\ge 0$), where $q_1, \ldots, q_\tau$ is a fixed set of
coprime integers. Clearly, such sequences grow subexponentially;
Tijdeman \cite{tij} proved that
\begin{equation}\label{subex}
n_{k+1}-n_k\ge \frac{n_k}{(\log n_k)^\alpha}
\end{equation}
for some $\alpha>0$, i.e. the growth of $(n_k)_{k \geq 1}$ is
almost exponential. Hardy-Littlewood-P\'olya sequences have
remarkable probabilistic and ergodic properties. Nair \cite{na}
proved that if $f$ is 1-periodic and integrable in $(0, 1)$, then
$$
\lim_{N\to\infty} \frac{1}{N} \sum_{k=1}^N f(n_kx)=\int_0^1 f(t)
dt\quad \textup{a.e.}
$$
Philipp \cite{ph1994} showed that the discrepancy of $\{n_kx\}$
satisfies the law of the iterated logarithm
\begin{equation}\label {hplil}
\frac{1}{4\sqrt{2}} \le \limsup_{N \to \infty} \frac{N D_N(n_k x)}
{\sqrt{2 N \log \log N}}\le C \quad \textup{a.e.}
\end{equation}
where $C$ is a constant depending on the number of generators of
$(n_k)_{k \geq 1}$. Recently, Fukuyama  and Nakata \cite{ft3}
succeeded in computing the limsup in (\ref{hplil}). Fukuyama and
Petit \cite{fupe} also showed that the central limit theorem
$N^{-1/2} \sum_{k=1}^N f(n_kx)\to_d N(0, \gamma_f^*)$ holds with
\begin{equation}\label{funaconst}
\gamma^*_f=\sum_{k, l: (n_k, n_l)=1}\int_0^1 f(n_kx)f(n_lx)dx.\\
\end{equation}
The Diophantine properties of $(n_k)_{k \geq 1}$ have been studied
in great detail in recent years; Amoroso and Viada \cite{am}
showed recently that Hardy-Littlewood-P\'olya sequences satisfy
condition ${\mathbf A}_p$ for any $p\ge 2$ with $C_p=\exp(p^6)$.
This is a very deep result, involving a substantial sharpening of
the subspace theorem of Evertse, Schlickewei and Schmidt (see
\cite{ev}). Again, the limit $\gamma$ in (\ref{gamma}) depends on
the permutation $\sigma$.\\

Since verifying condition $\mathbf{A}_p$ for a concrete
subexponential sequence $(n_k)_{k \geq 1}$ is difficult, it is
worth looking for Diophantine conditions which are strong enough
to imply the permutation-invariant CLT and LIL, but which hold for
a sufficiently large class of subexponential sequences.
Such a Diophantine condition $\boldsymbol{A_\omega}$ will be given
below. Actually, we will see that in a certain statistical sense,
$\boldsymbol{A_\omega}$ is satisfied for ``almost all'' sequences
$(n_k)_{k \geq 1}$ and thus the permutation-invariant CLT and LIL
are the ``typical'' behavior of sequences $f(n_kx)$. Given a
nondecreasing sequence $\boldsymbol{\omega}=(\omega_1, \omega_2,
\ldots)$ of positive integers tending to $+\infty$, we say that an
increasing sequence $(n_k)_{k \geq 1}$ of positive integers
satisfies

\bigskip\noindent
{\bf Condition $\boldsymbol{A_\omega}$}, if the Diophantine
equation
\begin{equation*}
a_1 n_{k_1} +\ldots +a_p n_{k_p}=0, \qquad k_1 < \ldots <k_p, \ \, 2
\leq p \leq \omega_N, \ \, |a_1|, \ldots, |a_p|\le N^{\omega_N}
\end{equation*}
has no nondegenerate solutions with $k_p>N$ (degenerate solutions
are
solutions where proper subsums vanish). \\

\begin{theorem}\label{th1a} Let {\boldmath{$\omega$}}=$(\omega_1,
\omega_2, \ldots)$ be a nondecreasing sequence tending to
$+\infty$ and let $(n_k)$ be an increasing sequence of positive
integers satisfying condition $\boldsymbol{A_\omega}$. Let $f$
satisfy (\ref{fcond}), let $\sigma : {\mathbb N}\to{\mathbb N}$ be
a permutation of the positive integers and assume that
\begin{equation}\label{dn}
d_N^2:=\int_0^1 \left(\sum_{k=1}^N f(n_{\sigma(k)} x)\right)^2\,
dx\ge CN \qquad (N\ge N_0)
\end{equation} for some constant $C>0$.
Then we have
\begin{equation}\label{permclt}
d_N^{-1}\sum_{k=1}^N f(n_{\sigma(k)} x) \stackrel{\mathcal D}
{\rightarrow} \mathcal{N}(0, 1).
\end{equation}
If $\omega_N \geq (\log N)^{\alpha}$ for some $\alpha>0$, then we
also have
\begin{equation} \label{permlildn}
\limsup_{N \to \infty} \frac{1}{\left( 2 d_N^2 \log \log
d_N^2\right)^{1/2}}
 \sum_{k=1}^N f(n_{\sigma(k)} x) = 1 \quad \textup{a.e.}
\end{equation}
\end{theorem}

Fix $\omega_N\to\infty$. We show that,
in a certain statistical sense, ``almost all'' sequences $n_k\le
k^{\omega_k}$ satisfy condition $\boldsymbol{A_\omega}$. To make
this precise, we need to define a probability measure over the set
of such sequences, or, equivalently, a natural random procedure to
generate such sequences. Clearly, the simplest procedure is to
choose $n_k$ independently and uniformly from the integers in the
interval $I_k=[1, k^{\omega_k}]$ $(k=1, 2, \ldots)$. Denote the so
obtained measure by $\mu$.

\begin{theorem}\label{aomega0}
With probability one with respect to $\mu$, the sequence $(n_k)_{k
\geq 1}$ satisfies condition $\boldsymbol{A_\omega}$.
\end{theorem}

As an immediate consequence, we get

\begin{theorem}\label{aomega}
With probability 1 with respect to $\mu$, the sequence $(n_k)_{k
\geq 1}$ obeys the central limit theorem (\ref{permclt}) with
$d_N = \|f\| \sqrt{N}$, and if $\omega_N\ge (\log N)^\alpha$ for
some $\alpha>0$, $(n_k)$ also satisfies the law of the iterated
logarithm (\ref{permlildn}) with $d_N = \|f\| \sqrt{N}$.
\end{theorem}

Clearly, for slowly increasing $(\omega_k)$ the so obtained
sequence $(n_k)$ grows almost polynomially (as a comparison,
Hardy-Littlewood-P\'olya sequences grow almost exponentially by
(\ref{subex})). We do not know if there exist polynomially growing
sequences $(n_k)_{k \geq 1}$ satisfying the permutation-invariant
CLT and LIL. As a simple combinatorial argument shows, sequences
$(n_k)_{k \geq 1}$ satisfying $\mathbf{A}_p$ for all $p\ge 2$
cannot grow polynomially.

\end{document}